\documentclass[10pt,reqno]{amsart}
\usepackage{amsmath}
\usepackage{amsfonts}
\usepackage{amssymb}
\usepackage{amsthm}
\usepackage{amscd}

\usepackage{a4wide}

\theoremstyle{remark}

\theoremstyle{definition}

\begin{document}
\large
\title{Proof of Bertrand's Postulate for $n\geq 6$}
\author{Manoj Verma}
\address{MANOJ VERMA: School of Physical Sciences, Jawaharlal Nehru University, New Delhi 110067 India}
\email{m.infinite@gmail.com}

\begin{abstract} We add a few ideas to Erd\H{o}s's proof of Bertrand's Postulate to produce one using a little calculus but requiring direct check only for $n\leq 5$ and one without using calculus and requiring direct check only for $n\leq 12$. The proofs can be presented to high school students.\\
\\
\end{abstract}
\maketitle
\noindent In 1845 J. F. L. Bertrand conjectured that for any real number $x>1$ there is a prime number $p$ satisfying $x<p<2x.$ It is equivalent to the same statement restricted to natural numbers $n>1$ or to primes. (Given $x\geq 2$, let $q$ be the largest prime less than or equal to $x$ and let $p$ be the prime satisfying $q<p<2q.$ Then $p>x$ by maximality of $q$ and $p<2q\leq 2x.$ For $1<x<2$, take $p=2.$) It was first proved by Chebyshev in 1852. In 1932 Paul Erd\H{o}s \cite{Erdos} gave a proof that does not require calculus. He showed that if there is no prime between $n$ and $2n$ then we must have $2^{2n/3}<(2n)^{1+\sqrt{2n}}$ and proved that this cannot hold for $n\geq 4000.$ Indeed, this inequality holds for $n=467 $ but does not hold for $n\geq 468.$ Bertrand's postulate can be checked to be true for $n\leq 3999$ using the primes 2, 3, 5, 7, 13, 23, 43, 83, 163, 317, 631, 1259, 2509 and 4001, each prime in the sequence being less than twice the previous one. The proof of Bertrand's postulate given in the widely used Number Theory textbook by Niven, Zuckerman and Montgomery \cite{NZM} uses calculus and requires direct check for $n\leq 1599.$ By adding a few simple observations to Erd\H{o}s's proof, we give a proof using a little calculus that requires direct check only for $n\leq 5$ and knowing that 2, 3, 5 and 7 are primes while 4, 6 and 8 are not  and no knowledge of primality of any other number. We also give a proof without calculus that requires direct check only for $n\leq 12$ and knowing that 2, 3, 5, 7 and 13 are primes while 4, 6, and 8 are not and no knowledge of primality of any other number.\\
\indent For a real number $x$, define $\pi(x)$ to be the number of primes less than or equal to $x$ and denote by $[x]$ the greatest integer less than or equal to $x$.\\\\
\noindent {\bf Lemma 1} We have the following bounds for $\binom{2n}{n}=\frac{(2n)!}{n!n!}$:
$$\frac{4^n}{\sqrt{3n+1}}\geq \frac{4^n}{\sqrt{\frac{28}{9}n+\frac{8}{9}}}\geq \binom{2n}{n}=\frac{(2n)!}{n!n!} \geq \frac{4^n}{\sqrt{\frac{16}{5}n+\frac{4}{5}}} \geq \frac{4^n}{\sqrt{4n}} \mbox{ with equality iff } n=1.$$
\noindent {\bf Proof} Consider the expression $\frac{4^n}{\sqrt{xn+y}}.$ It equals $\binom{2n}{n}$ when $n=1$ if $x+y=4.$ For $n\geq 1,$ on increasing $n$ by 1 it gets multiplied by $\frac{4\sqrt{xn+y}}{\sqrt{xn+x+y}}$ while $\binom{2n}{n}$ gets multiplied by $\frac{(2n+1)(2n+2)}{(n+1)^2}=\frac{2(2n+1)}{n+1}.$ Now
$$\frac{4\sqrt{xn+y}}{\sqrt{xn+x+y}}>\mbox{ or }=\mbox{ or }<\frac{2(2n+1)}{n+1}$$
according as 
$$4(xn+y)(n+1)^2>\mbox{ or }=\mbox{ or }<(2n+1)^2(xn+x+y).$$
The left hand side equals $4xn^3+(8x+4y)n^2+(4x+8y)n+4y$ while the right hand side equals $4xn^3+(8x+4y)n^2+(5x+4y)n+x+y$. The left hand side is greater than the right hand side for $n\geq 1$ provided $(4y-x)n+3y-x=(4y-x)(n-1)+7y-2x\geq 0$ for $n\geq 1$ (so we can take $x=28/9, y=8/9$ or $x=3, y=1$). The left hand side is less than the right hand side for $n\geq 1$ provided $(4y-x)n+3y-x=(4y-x)(n-1)+7y-2x\leq 0$ for $n\geq 1$ (so we can take $x=32/5, y=8/5$ or $x=4, y=0$). This completes the proof.\\\\
{\bf Lemma 2} We have
$$\prod_{p\leq x}p< \frac{4^x}{6x}\mbox{ for any real number } x\geq 4.$$
\noindent {\bf Proof} For $4\leq x<5$, we have $\prod_{p\leq x}p=6<8=\frac{4^4}{32}< \frac{4^x}{6x+2}.$ Note that for any real number $x\geq 1$ we have
 $$\frac{4^{x+\frac{1}{3}}}{6(x+\frac{1}{3})+2}=\frac{4^{x}}{6x+2}\frac{4^{1/3}}{(6x+4)/(6x+2)}>\frac{4^{x}}{6x+2}$$
since $\frac{6x+4}{6x+2}=1+\frac{1}{3x+1}\leq \frac{5}{4}=\sqrt[3]{\frac{125}{64}}<\sqrt[3]{4}.$ It follows that for any natural number $n$
 $$\frac{4^{n}}{6n+2}<\frac{4^{n+\frac{1}{3}}}{6(n+\frac{1}{3})+2}< \frac{4^{n+\frac{2}{3}}}{6(n+\frac{2}{3})+2}<\frac{4^{n+1}}{6(n+1)+2}.$$
If we prove that $\prod_{p\leq n}p< \frac{4^n}{6n+2}\mbox{ for all natural numbers } n\geq 5,$ it would follow that
$$\prod_{p\leq x}p=\prod_{p\leq n}p<\frac{4^{n}}{6n+2}<\frac{4^x}{6x} \mbox{ for } n\geq 5 \mbox{ and } n\leq x\leq n+\frac{1}{3},$$
$$\prod_{p\leq x}p=\prod_{p\leq n}p<\frac{4^{n}}{6n+2}<\frac{4^{n+\frac{1}{3}}}{6n+4}<\frac{4^x}{6x} \mbox{ for } n\geq 5 \mbox{ and } n+\frac{1}{3}\leq x\leq n+\frac{2}{3},$$
and
$$\prod_{p\leq x}p=\prod_{p\leq n}p<\frac{4^{n}}{6n+2}<\frac{4^{n+\frac{2}{3}}}{6n+6}<\frac{4^x}{6x} \mbox{ for } n\geq 5 \mbox{ and } n+\frac{2}{3}\leq x\leq n+1$$
completing the proof of the proposition.\\
(If we use calculus, we just note that the right-hand-side of the inequality is an increasing function of $x$ for $x\geq 1$ as its derivative $(4^x/6x)(\log 4 -1/x)$ is positive for $x\geq 1$. Thus it is enough to prove that
$$\prod_{p\leq n} p < \frac{4^n}{6n}$$
for any natural number $n\geq 4$ as then it will follow for any real number $x\geq 4$ that
$$\left.\prod_{p\leq x} p=\prod_{p\leq [x]} p < \frac{4^{[x]}}{6[x]}\leq \frac{4^x}{6x}. \right)$$
\indent We shall use induction in the following form: Let $n_0\geq 2$ be a natural number. If a statement $P(n)$ about natural numbers $n$ is true for $n_0\leq n\leq 2n_0-2$ and $P(2n-1)$ and $P(2n)$ are true whenever $n\geq n_0$ and $P(n)$ is true, then $P(n)$ is true for all $n\geq n_0$.\\
\indent The inequality $\prod_{p\leq n}p< \frac{4^n}{6n+2}$  can be checked directly for $5\leq n\leq 8$: $$\prod_{p\leq 5}p=\prod_{p\leq 6}p=30< 32=\frac{4^5}{32}<\frac{4^6}{38}$$ and $$\prod_{p\leq 7}p=\prod_{p\leq 8}p=210<\frac{4096}{11}=\frac{4^7}{44}<\frac{4^8}{50}.$$ Assuming it true for $n=m\geq 5$ we shall prove that it holds for $n=2m-1$ and $n=2m$ as well, completing the proof by induction. Since $$\binom{2m-1}{m-1}=\frac{(m+1)(m+2)\ldots (2m-1)}{1\cdot 2\cdot \ldots \cdot (m-1)}=\frac{1}{2}\binom{2m}{m}$$ is a positive integer divisible by any prime $p$ with $m+1\leq p\leq 2m-1$, it is greater than or equal to the product of primes $p$ with $m+1\leq p\leq 2m-1.$ Hence, for $m\geq 5,$ we have
$$\prod_{m<p\leq 2m-1}p=\prod_{m<p\leq 2m}p\leq \binom{2m-1}{m-1}=\frac{1}{2}\binom{2m}{m}<\frac{4^{m}}{2\sqrt{3m+1}}\leq \frac{4^{m-1}}{2}$$
and thus using the induction hypothesis $$\prod_{p\leq 2m-1}p=\prod_{p\leq 2m}p<\frac{4^m}{6m+2}\cdot\frac{4^{m-1}}{2}=\frac{4^{2m-1}}{12m+4}<\min\left(\frac{4^{2m-1}}{6(2m-1)+2}, \frac{4^{2m}}{6(2m)+2}\right).$$
This completes the proof.\\\\
\indent We recall Legendre's formula for the highest power of a prime $p$ dividing a factorial $n!$:
$$n!=\prod_{p\leq n}p^{a_{p}} \mbox{ where } a_{p}=\sum_{i=1}^{\infty}[n/p^{i}]=\sum_{i=1}^{r_p}[n/p^{i}]$$
where $r_p$ is the largest integer $r$ such that $p^{r}\leq n$. 
We rephrase Bertrand's postulate.\\\\
\noindent {\bf Bertrand Postulate} For any natural number $m\geq 2$, $$\prod_{m<p\leq 2m} p >1.$$
\noindent {\bf Proof}  From Legendre's formula for the exponent of the highest power of a prime $p$ dividing a factorial and Lemma 1 we have
$$
\frac{2^{2m-1}}{\sqrt{m}}\leq \frac{(2m)!}{m!m!}=\prod_{p\leq 2m}p^{\left(\sum_{i=1}^{r_p}\left[\frac{2m}{p^i}\right]-2\left[\frac{m}{p^i}\right]\right)}.
$$
where $r_p$ is the largest integer $r$ such that $p^{r}\leq 2m.$
For any real number $x$ we have $[2x]-2[x]=0$ or 1. We note that for $m+1\leq p\leq 2m$ there is only one term in the sum and it equals 1. For $m\geq 5$, $(2m/3)^2>2m$ so there is only one term in the sum; for $2m/3<p\leq m$, it is $2-2= 0$, for $\sqrt{2m}<p\leq 2m/3$ it is at most 1. For $p\leq \sqrt{2m}$ the exponent is at most $r_p$ and the corresponding factor is still at most $2m.$ Hence, for $m\geq 5$, we have
$$\frac{2^{2m-1}}{\sqrt{m}}\leq \frac{(2m)!}{m!m!}\leq \left(\prod_{m<p\leq 2m}p\right)\cdot \left(\left(\prod_{p\leq 2m/3}p\right)\left(\prod_{p\leq \sqrt{2m}}p\right)^{-1}\right)\cdot (2m)^{\pi(\sqrt{2m})}.$$
Thus
$$\prod_{m<p\leq 2m}p\geq \frac{2^{2m-1}}{\sqrt{m}}\left(\prod_{p\leq 2m/3}p\right)^{-1}\left(\prod_{p\leq \sqrt{2m}}p\right)\cdot (2m)^{-\pi(\sqrt{2m})}.$$
Now $\prod_{p\leq 2m/3}p<\frac{2^{4m/3}}{4m}$ by Proposition 1 provided $2m/3\geq 4$, i.e., $m\geq 6$ and $\prod_{p\leq \sqrt{2m}}p\geq 2^{\pi(\sqrt{2m})}$. Thus, for $m\geq 6,$ we have
$$\prod_{m<p\leq 2m}p> 2^{(2m/3)+1}\sqrt{m}\cdot m^{-\pi(\sqrt{2m})}.$$
Since 1 is not a prime and all the primes except 2 or odd, $\pi(x)\leq$ (number of odd natural numbers $\leq x)\leq$ (number of even natural numbers $\leq x+1)=[(x+1)/2]\leq (x+1)/2.$ Thus
\begin{equation}  \label{e3}
\prod_{m<p\leq 2m}p> 2^{(2m+3)/3}\cdot \sqrt{m}\cdot m^{-({\sqrt{2m}+1})/2}
= 2^{(2m+3)/3}\cdot m^{-\sqrt{2m}/2}.
\end{equation}
Thus $\prod_{m<p\leq 2m}p>1$ provided 
$m^3\leq 2^{\frac{(4m+6)}{\sqrt{2m}}}=2^{2\sqrt{2m}+\frac{6}{\sqrt{2m}}}$, i.e.,
\begin{equation} \label{e1}
 (\sqrt{2m})^6<2^{2\sqrt{2m}+\frac{6}{\sqrt{2m}}+3}.
\end{equation}
This inequality holds for all $m\geq 1$ (Lemma 3 below) but we assumed $m\geq 6$ in an earlier step.
For $m\leq 5$ we can check the postulate directly by using the primes 2, 3, 5 and 7. For $m=2$ take $p=3$; for $m=3, 4$ take $p=5$; for $m=5$ take $p=7.$ ($p=7$ works for $m=6$ as well.)\\\\
\noindent {\bf Lemma 3} For any real number $x>0$ we have
\begin{equation} \label{e2}
x^6<2^{2x+3+\frac{6}{x}}.
\end{equation}
\noindent {\bf Proof 1 (Using calculus)} The inequality is equivalent to $2x+3+\frac{6}{x}-\frac{6\log x}{\log 2}>0$ for $x>0.$ The function $f_{\alpha}(x)=\alpha x+3+\frac{6}{x}-\frac{6\log x}{\log 2}$ has the derivative $\alpha-\frac{6}{x^2}-\frac{6}{x\log 2}$ which increases steadily with $x$ and becomes zero at $x=x_0$ where $x_0$ is the solution of the equation $\alpha=\frac{6}{x_0^2}+\frac{6}{x_0\log 2}$ and thus $f_{\alpha}$ attains at $x=x_0$ its minimum value of $$
\left(\frac{6}{x_0^2}+\frac{6}{x_0\log 2}\right)x_0+3+\frac{6}{x_0}-\frac{6\log x_0}{\log 2}=\frac{12}{x_0}+\frac{6}{\log 2}+3-\frac{6\log x_0}{\log 2}.$$
Since $2x+3+\frac{6}{x}-\frac{6\log x}{\log 2}=(2-\alpha)x+f_{\alpha}(x)$, the inequality required would follow if we could choose $x_0$ such that $\alpha\leq 2$ and the minimum value attained by $f_{\alpha}$ is positive. Taking $x_0=5$ we have $\alpha=\frac{6}{25}+\frac{6}{5\log 2}=1.971234...<2$ and the minimum value of $f_{\alpha}$ is $f_{\alpha}(5)=\frac{12}{5}+\frac{6}{\log 2}+3-\frac{6\log 5}{\log 2}=0.1246...>0.$  This completes the proof of Bertrand's postulate using calculus.\\\\
\indent One salient feature of Erd\H{o}s's original proof of Bertrand's postulate is that it does not use caluculus. We can prove the inequality (\ref{e2}) without using calculus for all $x>0$ if we are willing to do some computation to verify the inequality for small values of $x$ as in Proof 4 below. Alternatively, we note that the inequality (\ref{e2}) for $x\geq x_0$ implies the inequality (\ref{e1}) for $m\geq x_0^{2}/2$ and we can check (\ref{e1}) or Bertrand's postulate directly for the smaller values of $m.$\\
\indent We shall use the following variation of induction: Let $\alpha$ be a positive real number. If a statement $P(x)$ about real numbers $x$ is true for $x_0\leq x< x_0+\alpha$ and $P(x+\alpha)$ is true whenever $x\geq x_0$ and $P(x)$ is true, then $P(x)$ is true for all real numbers $x\geq x_0.$\\
\indent We observe that for $x\geq 2$, $2x\geq4>3\geq 6/x$, $2x$ increases while $6/x$ decreases as $x$ increases hence $2x-(6/x)+3=(\sqrt{2x}-\sqrt{6/x})^2+2\sqrt{12}+3$ increases with $x$.\\
\indent For $n\geq 2$, we observe that
$$\left(1+\frac{1}{n}\right)^n>\left(1+\frac{1}{n}\right)\cdot \left(1+\frac{1}{n+1}\right)\cdots \left(1+\frac{1}{2n-1}\right)$$$$=\frac{n+1}{n}\cdot \frac{n+2}{n+1}\cdots \frac{2n}{2n-1}=2$$
and
$$\left(1+\frac{1}{2n-1}\right)^n<\left(1+\frac{1}{2n-1}\right)\cdot \left(1+\frac{1}{2n-2}\right)\cdots \left(1+\frac{1}{n}\right)$$$$=\frac{2n}{2n-1}\cdot \frac{2n-1}{2n-2}\cdots \frac{n+1}{n}=2.$$ Hence
$$1+\frac{1}{2n}< 1+\frac{1}{2n-1}\leq 2^{\frac{1}{n}}\leq 1+\frac{1}{n} \mbox{ with equality iff } n=1.$$\\
\noindent {\bf Proof 2 (Without calculus)}  For $5\leq x\leq 5.1$ we have $x^6\leq (5.1)^6=((5.1)^3)^2<(2^{7.1})^2=2^{14.2}=2^{2(5)+3+\frac{6}{5}}\leq 2^{2x+3+\frac{6}{x}}$ since $2^{7.1}=2^{7}\cdot 2^{\frac{1}{10}}>128(1+\frac{1}{20})>128+6=134$ while $5.1^3=5^3+3\cdot 5^2\cdot (0.1)+3\cdot 5\cdot (0.01)+0.001<125+8=133.$ For $x\geq 5$, on increasing $x$ by 0.1, $x^6$ increases by a factor of
$$\frac{(x+0.1)^6}{x^6}=\left(1+\frac{0.1}{x}\right)^6\leq \left(1+\frac{1}{50}\right)^6$$
while $2^{2x+3+\frac{6}{x}}$ increases by a factor of
$$\frac{2^{2(x+0.1)+3+\frac{6}{x+0.1}}}{2^{2x+3+\frac{6}{x}}}=2^{0.2-\frac{0.6}{x(x+0.1)}}\geq 2^{0.2-\frac{0.6}{5(5.1)}}= 2^{\frac{1}{5}-\frac{2}{85}}
=2^{\frac{3}{17}}=8^{\frac{1}{17}}.$$
As $\left(\left(1+\frac{1}{50}\right)^6\right)^{17}=\left(1+\frac{1}{50}\right)^2\cdot\left(1+\frac{1}{50}\right)^{50\cdot 2}<\left(1+\frac{1}{44}\right)^2\cdot e^2<\left(\frac{45}{44}\right)^2\cdot \left(\frac{11}{4}\right)^2=\frac{2025}{256}<8$, the inequality (\ref{e2}) holds for $x\geq 5$ by induction and  the inequality (\ref{e1}) holds for $m\geq 5^2/2=12.5$,  i.e. for $m\geq 13.$ For $7\leq m\leq 12,$ take $p=13$. Here, we have used the following. For any real number $x>0$, $e^{x}=1+x+\frac{x^2}{2!}+\ldots>1+x$ hence $(1+x)^{\frac{1}{x}}<e.$ In particular, $\left(1+\frac{1}{n}\right)^n<e$ for any natural number $n$. Also, $e=1+\frac{1}{1!}+\frac{1}{2!}+\frac{1}{3!}+\ldots<1+\frac{1}{1!}+\frac{1}{2!}\left(1+\frac{1}{3}+\frac{1}{3^2}+\frac{1}{3^3}+\ldots\right)=1+1+\frac{1}{2}\cdot\frac{3}{2}=\frac{11}{4}.$\\\\
\noindent {\bf Proof 3 (Without calculus or the infinite series for $e^x$)} For $5.8\leq x\leq 6$ we have $x^6\leq 6^6<(2^{2.6})^6=2^{15.6}=2^{2(5.8)+3+\frac{6}{6}}\leq 2^{2(5.8)+3+\frac{6}{5.8}}\leq 2^{2x+3+\frac{6}{x}}$ as $(2^{2.6})^5=2^{13}=8192>6^5=7776$. For $x\geq 5.8$, on increasing $x$ by 0.2, $x^6$ increases by a factor of
$$\frac{(x+0.2)^6}{x^6}=\left(1+\frac{0.2}{x}\right)^6\leq \left(1+\frac{1}{29}\right)^6<\frac{30}{29}\cdot \frac{29}{28}\cdot \frac{28}{27}\cdot \frac{27}{26}\cdot \frac{26}{25}\cdot \frac{25}{24}= \frac{30}{24}=\frac{5}{4}$$
while $2^{2x+3+\frac{6}{x}}$  increases by a factor of 
$$\frac{2^{2(x+0.2)+3+\frac{6}{x+0.2}}}{2^{2x+3+\frac{6}{x}}}=2^{0.4-\frac{1.2}{x(x+0.2)}}\geq 2^{0.4-\frac{1.2}{(5.8)(6.0)}}= 2^{\frac{2}{5}-\frac{1}{29}}.$$
Now $(\frac{5}{4})^5=\frac{3125}{1024}=3\frac{53}{1024}<3\frac{1}{3}$ while $(2^{\frac{2}{5}-\frac{1}{29}})^5=4\cdot (2^{-\frac{5}{29}})>4\cdot 2^{-\frac{1}{5}} >4\cdot (1+\frac{1}{5})^{-1}=4\cdot \frac{5}{6}=3\frac{1}{3}.$ Hence the inequality (\ref{e2}) holds for $x\geq 5.8$ by induction and  the inequality (\ref{e3}) holds for $m\geq 5.8^2/2=16.82$,  i.e. for $m\geq 17.$ For 
$7\leq m\leq 12$ take $p=13$; for $13\leq m\leq 16$ take $p=17.$\\\\
\noindent {\bf Proof 4 (Without calculus but with some calculation)} For $0<x\leq 2$, $x^6\leq 2^6\leq2^{\frac{6}{x}+3}<2^{2x+\frac{6}{x}+3}.$ For $2\leq x\leq 3$, $x^6\leq 3^6=729<1024=2^{10}=2^{2(2)+\frac{6}{2}+3}\leq 2^{2x+\frac{6}{x}+3}.$  For $3\leq x\leq 3.5$, $x^6\leq (3.5)^6=(\frac{343}{8})^2<43^2=1849<2048=2^{11}=2^{2(3)+\frac{6}{3}+3}\leq 2^{2x+\frac{6}{x}+3}.$ For $3.5\leq x\leq 3.75$, $x^6\leq (3.75)^6=(\frac{3375}{64})^2<53^2=2809<2828=(1.4)(2)(1010)<2^{0.5+1+10}<2^{82/7}=2^{2(3.5)+\frac{6}{3.5}+3}\leq2^{2x+\frac{6}{x}+3}.$ For $3.75\leq x\leq 4$, $x^6\leq 4^6=2^{12}<2^{12.1}=2^{2(3.75)+\frac{6}{3.75}+3}<2^{2x+\frac{6}{x}+3}.$ In the same way we can verify the inequality for $4\leq x\leq 6$ by partitioning the interval $4\leq x\leq 6$ using the points 4, 4.2, 4.4, 4.55, 4.7, 4.85, 5, 5.15, 5.3, 5.45, 5.6, 5.8 and 6 using a calculator. For $x\geq 4.5$  on increasing $x$ by 1.5, $x^6$ increases by a factor of
$$\frac{(x+1.5)^6}{x^6}=\left(1+\frac{1.5}{x}\right)^6\leq \left(1+\frac{1.5}{4.5}\right)^6=\left(\frac{4}{3}\right)^6=\frac{4096}{729}<6.$$
while $2^{2x+3+\frac{6}{x}}$ increases by a factor of
$$\frac{2^{2(x+1.5)+\frac{6}{x+1.5}+3}}{2^{2x+\frac{6}{x}+3}}=2^{3-\frac{9}{x(x+1.5)}}\geq 2^{3-\frac{9}{(4.5)(6)}}
=2^{\frac{8}{3}}=\sqrt[3]{256}>\sqrt[3]{216}=6.$$
Hence the inequality holds for $x\geq 6$ as well by induction.\\\\
{\bf Remark 1} From $(\ref{e3})$, for $m\geq 6$, $\prod_{m<p\leq 2m}p> 2m$ 
provided $2^{(2m+3)/3}\cdot m^{-\sqrt{2m}/2}\geq 2m$, i.e., $\frac{2m\log 2}{3}-\left(1+\frac{\sqrt{2m}}{2}\right)\log m\geq 0.$ The function $f(x)=\frac{2x\log 2}{3}-\left(1+\frac{\sqrt{2x}}{2}\right)\log x$ is positive at $x=53$ and its derivative $\frac{2\log 2}{3}-\frac{\log x}{2\sqrt{2x}}-\frac{1}{x}-\frac{1}{\sqrt{2x}}$ increases steadily for $x\geq e^2=7.389...$ and is positive at $x=19.$ By a direct check for $x<53$ (for $x$ in the interval $[2.5, 3), [3.5, 5), [5.5, 7), [7, 11), [11, 17), [17, 29), [29, 47), [47, 83)$ respectively, take the prime pair to be $\{3, 5\}, \{5, 7\}, \{7, 11\}, \{11, 13\}, \{17, 19\}, \{29, 31\}, \{47, 53\}, \{83, 89\}$\\ respectively), we see that for any real number $x\in [2.5, 3)\cup [3.5, 5) \cup [5.5, \infty)$ there are at least two primes $p$ satisfying $x<p\leq 2x$. For $m=4$ and $m\geq 6$ there are at least two primes $p$ satisfying $m<p<2m$ and hence at least one prime $p$ satisfying $m<p\leq 2m-3$.\\\\
{\bf Remark 2} From $(\ref{e3})$, the number of primes $p$ satisfying $m<p\leq 2m$ is $$> \frac{(2m+3)\log 2}{3\log 2m}-\frac{\sqrt{2m}\log m}{2\log 2m}>\frac{2m\log 2}{3\log 2m}-\frac{\sqrt{2m}\log m}{2\log 2m}=\frac{2m\log 2}{3 \log 2m}\left(1-\frac{3\log m}{2\sqrt{2m}\log 2}\right)$$ and thus at least $\frac{0.3 m}{\log 2m}$ for $m\geq 1000$ and at least $\frac{0.4m}{\log 2m}$ for $m\geq 11292.$

\end{document}